\documentclass[11pt]{article}

\usepackage{mathptmx,eucal,amsmath,amscd,amssymb,amsthm,xspace}
\usepackage{hyperref,amssymb,color, url,fancyhdr}

\theoremstyle{plain}

\newtheorem{thm}{Theorem}[section]

\newtheorem{cor}[thm]{Corollary}

\newtheorem{conj}{Conjecture}

\newtheorem{lemma1}[thm]{Lemma}
\newtheorem{lemma2}[thm]{Lemma}
\newtheorem {lemma3} [thm] {Lemma}

\theoremstyle{definition}

\newtheorem*{q1}{Question 1}
\newtheorem*{q2}{Question 2}

\def\hat{\widehat}

\fancypagestyle{plain}

\begin{document}
\bibliographystyle{plain}


\title{\textbf{Simple closed curves, word length, and nilpotent quotients of free groups}}
\author{Khalid Bou-Rabee and Asaf Hadari}
\maketitle


\begin{abstract}
We consider the fundamental group $\pi$ of a surface of finite type equipped with the infinite generating set consisting of all simple closed curves. We show that every nilpotent quotient of $\pi$ has finite diameter with respect to the word metric given by this set. This is in contrast with a result of Danny Calegari that shows that $\pi$ has infinite diameter with respect to this set. Furthermore, we give a general criterion for a finitely generated group equipped with a generating set to have this property.
\end{abstract}

\smallskip\smallskip


\section{Introduction}

A surface of finite type is a surface whose fundamental group is finitely generated. Given such a surface there is no canonical choice of generating set. If one wishes to define a suitably canonical generating set of a geometric nature then it becomes necessary to consider infinite generating sets. One such set is the set of all elements whose conjugacy class can be represented by a simple closed curve. These are in some sense the simplest elements of the fundamental group, and are thus a natural choice for a generating set.

Benson Farb posed the question whether the fundamental group, endowed with the word metric given by this set, has finite diameter. This question was answered negatively by Danny Calegari \cite{C08}. In this paper our goal is to investigate the same question for some quotients of the fundamental group. In contrast with Calegari's result, we find the following.

\begin{thm} \label{thm1.1}
Let $\Sigma$ be a surface of finite type, $\pi = \pi_{1}(\Sigma)$, and $\mathcal{S} \subset \pi$ be any generating set containing at least one element in each conjugacy class that is represented by a nonseparating simple closed curve. Let $\rho: \pi \to N$ be a homomorphism into any nilpotent group. Then $\rho(\pi)$ has finite diameter in the word metric with respect to the set $\rho(\mathcal{S})$.
\end{thm}

Note that in surfaces of genus $> 1$, $\pi$ has many nilpotent quotients of every degree of nilpotency. Furthermore, it is residually nilpotent, that is for every $x \in \pi$ there is some nilpotent quotient $q: \pi \to N$ such that $q(x) \neq 1$.

We say that a group $G$ is \textit{nilpotent-bounded with respect to the set S} if any nilpotent quotient of $G$ has finite diameter with respect to the word metric given by the image of $S$. As part of the proof we prove the following more general result.

\begin{thm} \label{thm1.2}
Let $G$ be a finitely generated group, and let $S\subset G$ be a generating set such that $G/[G,G]$ has finite diameter with respect to the word metric given by $S$. Then $G$ is nilpotent-bounded with respect to $S$.
\end{thm}

Using Theorem \ref{thm1.2}, it is possible to find smaller generating sets for which $\pi$ is nilpotent bounded. We give one such set here, but it is relatively simple to find many of them. In order to do so, we need a simple corollary. 

\begin{cor}\label{cor1.3}
Let G be a finitely generate group. Let $H = H_{1}(G, \mathbb{Z}) \cong G / [G,G]$. Suppose that $H \cong H_{1} \oplus \ldots \oplus H_{k}$, and that for each $i = 1, \ldots, k$ we are given a set $S_{i} \subset \Sigma$ whose projection to $H$ is contained in $H_{i}$ and generates $H_{i}$ with finite diameter. Then $G$ is nilpotent bounded with respect to $S_{1} \cup \ldots \cup S_{k}$.
\end{cor}

An example of an application of Corollary \ref{cor1.3} is the following.  Let $\Sigma$ be an orientable  of genus $g>1$. It is common to choose a generating set for $\pi = \pi_{1}(\Sigma)$ of the form $S' = \{\alpha_{1}, \beta_{1}, \ldots, \alpha_{g}, \beta_{g}\}$ where all of the above are represented by simple closed curves, the geometric intersection number of $\alpha_{i}$ and $\beta_{i}$ is one and they can be realized disjointly from all the other curves. Let $\Gamma_{i} = <\alpha_{i}, \beta_{i}>$. The group $\Gamma_{i}$ is the fundamental group of an embedded torus. Let $H = H_{1}(\Sigma)$, and $H_{i}$ be the projection to $H$ of $\Gamma_{i}$. Then $H = H_{1} \oplus \ldots \oplus H_{g}$. Thus, if we let $\mathcal{S}$ be the any set containing at least one representative in each conjugacy class of a simple closed curves that lies in one of the $g$ embedded tori described above, then $\pi$ is nilpotent bounded with respect to $\mathcal{S}$.

\paragraph{Acknowledgements.}
The authors wish to thank their advisor - Benson Farb, for his interest, suggestions, comments.

\section{Nilpotent Groups and Lower Central Series.}

Given a group $\Gamma$, we define a decreasing sequence of subgroups of $\Gamma$ called \textit{the lower central series of } $\Gamma$ by the following rule:

$$\Gamma_{0} = \Gamma, \textrm{                                  } \Gamma_{n+1} = [\Gamma,\Gamma_{n}]. $$

\noindent
A group is nilpotent if $\Gamma_{n} = \left<1\right>$ for some $n$. A group is called $n$-step nilpotent if $L_{n} = 1$, and $L_{n-1} \neq 1$. For every $n$, the group $L_{n} := \Gamma/ \Gamma_{n}$ is a nilpotent group. These groups have the property that any nilpotent quotient of $G$ factors through one of the projections $\Gamma \to L_{n}$.

Let $A_{n} := \Gamma_{n} / \Gamma_{n+1}$. It is a standard fact that $A_{n} = Z(L_{n})$, the center of $L_{n}$. Furthermore, if $\Gamma$ is finitely generated then $A_{n}$ is also finitely generated. Given a generating set $S$ of $\Gamma$, the group $A_{n}$ is generated by the images of elements of the form $[a_{1}, \ldots, a_{n}]$ where  $a_{1}, \ldots, a_{n} \in S$ and $[a_{1}, \ldots, a_{n}]$ denotes a generalized commutator, i.e:

$$[a_{1}, \ldots, a_{n}] =  [\ldots [a_{1}, a_{2}], a_{3}],  \ldots, a_{n}]$$
\noindent
In the course of the proof, we require the following technical lemma about generalized commutators in nilpotent groups.

\begin{lemma3} \label{CommutatorLemma}
Let $\Gamma$ be any group, $n,k \in \mathbb{N}$, and $a_{1}, \ldots, a_{n} \in \Gamma$. Then:

$$[a_{1}, \ldots, a_{n}]^{k} \equiv_{n+1} ([a_{1}^{k}, \ldots, a_{n}])   $$

\noindent where $\equiv_{i}$ is understood as having equal images in $L_{i}$

\end{lemma3}

\paragraph{Proof.} First, recall that $A_{n} = Z(L_{n+1})$. Let $x \in \Gamma_{n - 1}$ and $y \in \Gamma$.
Note that $[x,y] \in \Gamma_{n}$. Thus we have that:

$$[x^{k},y] \equiv_{n+1} x^{k}yx^{-k}y^{-1} \equiv_{n+1} x^{k}y[x,y]^{k}y^{-1}x^{-k} \equiv_{n+1} [x,y]^{k}. $$

\noindent The last equality stems from the fact that $[x,y]^{k}$ is central in $L_{n+1}$ and thus is invariant under conjugation. Note that this proves the claim for the case $n=1$. We now proceed by induction.

\noindent
By the case $n = 1$ we have that:

$$[a_{1}, \ldots, a_{n}]^{k} \equiv_{n+1} [[a_{1}, \ldots, a_{n-1}], a_{n}]^{k} \equiv_{n+1} [[a_{1}, \ldots, a_{n-1}]^{k}, a_{n}].$$

\noindent
By induction, we can write:
$$[a_{1}, \ldots, a_{n-1}]^{k} \equiv_{n+1} [[a_{1}, \ldots, a_{n-2}]^{k}, a_{n-1}]\gamma_{n}, $$

\noindent
where $\gamma_{n} \in \Gamma_{n}$. Since the image of $\Gamma_{n}$ is central in $L_{n+1}$ we have that :
$$[[a_{1}, \ldots, a_{n-1}]^{k}\gamma_{n}^{-1}, a_{n}]  \equiv_{n+1} [a_{1}, \ldots, a_{n-1}]^{k}, a_{n}].$$
Proceeding similarly we get the claim of the lemma. $\Box$

\section{Proof of the Main Theorems} \label{AnotherSection}

\begin{lemma1} \label{SympLemma}
Let $n \in \mathbb{N}$ and let $e_{1}, \ldots, e_{2n}$ be the standard basis for $\mathbb{Z}^{2n}$. Then the set $\mathcal{S} = Sp_{2n}(\mathbb{Z})\cdot e_{1}$ generates $\mathbb{Z}^{2n}$ with finite diameter.
\end{lemma1}

\paragraph{Proof.} We prove this fact first for $n=1$. In this case $Sp_{2n}(\mathbb{Z}) = SL_{2}(\mathbb{Z})$. Given a vector
 $v = {\tiny \left(
                \begin{array} {cc}
                a   \\
                b   \end{array} \right)} \in \mathbb{Z}^{2}$ such that $\textrm{gcd}(a,b) = 1$, there exist $x,y \in \mathbb{Z}$ such that $ax +by = 1$. In this case

 $$A = \left(
                \begin{array}{cc}
                a & -y  \\
                b & x  \end{array} \right) \in SL_{2}(\mathbb{Z})$$ and $A \cdot e_{1} = v$ and thus $v \in \mathcal{S}$. For a general vector $v = \left(
                \begin{array} {cc}
                a   \\
                b   \end{array} \right)$ notice that

                $$v =  \left(
                \begin{array} {cc}
                a - 1   \\
                1   \end{array} \right) + \left(
                \begin{array} {cc}
                1   \\
                b - 1   \end{array} \right)$$ and that $\textrm{gcd}(1,a-1) = \textrm{gcd}(1,b-1) = 1$, and thus $v \in \mathcal{S} + \mathcal{S}$.

 Now consider the case $n > 1$. In this case we have that $D < Sp_{2n}(\mathbb{Z})$ where $D \cong \prod_{i = 1}^{n}SL_{2}(\mathbb{Z})$ is the group of matrices containing $n$ copies of $SL_{2}(\mathbb{Z})$ along the diagonal, and zeroes in all other entries. Notice further that $\hat{e} = e_{1} + e_{3} + \ldots + e_{2n -1} \in \mathcal{S}$. Given ${\tiny \left(
                \begin{array} {cc}
                a_{i}   \\
                b_{i}   \end{array} \right)_{i = 1}^{n}} \in \mathbb{Z}^{2n}$, by the $n = 1$ case there are $2n$ matrices $A_{1}, \ldots, A_{n}, B_{1}, \ldots B_{n} \in SL_{2}(\mathbb{Z})$ such that:

                $$A_{i}\cdot e_{1} = \left(
                \begin{array} {cc}
                a_{i} - 1   \\
                1   \end{array} \right) , \textrm{      } B_{i}\cdot e_{1} = \left(
                \begin{array} {cc}
                1   \\
                b_{i} - 1   \end{array} \right).$$

\noindent
Let $A = \textrm{diag}(A_{1}, \ldots, A_{n})$, $B= {diag}(B_{1}, \ldots, B_{n})$  and . Then $$ v = A \cdot \hat{e} + B \cdot \hat{e}.$$
\noindent
Thus $\mathbb{Z}^{2n}$ is generated by $\mathcal{S}$ with finite diameter.
$\Box$

\begin{lemma2} \label{BoundedGenerationLemma}
Let $\Gamma$ be a finitely generated group, and let $n \in \mathbb{N}$. Suppose that $\mathcal{S} \subset \Gamma$ generates $\Gamma$ and generates $L_{n}$ with finite diameter. Then $\mathcal{S}$ generates $L_{n+1}$ with finite diameter.

\end{lemma2}

\paragraph{Proof.}
By assumption, there exists a $N_{0}$ such that for any $w\in \Gamma$ there exist $s_{1}, \ldots s_{m} \in \mathcal{S}$ (with $m < N_{0}$) such that
$$(s_{1}\ldots s_{m})^{-1}w \in \Gamma_{n}.$$
\noindent
Thus, it is enough to show that the image of $\mathcal{S}$ in $L_{n+1}$  generates $A_{n}$ with finite diameter. The group $A_{n}$ is a finitely generated abelian group which is  generated by elements of the form $[s_{1}, \ldots, s_{n}]$ where  $s_{1}, \ldots s_{n} \in \mathcal{S}$. Choose $p$ such generators: $\gamma_{1}, \ldots, \gamma_{p}$.
\noindent
Consider $\gamma_{1} = [s_{1}, \ldots, s_{n}]$. Given any $k \in \mathbb{N}$, by Lemma \ref{CommutatorLemma}, we have that $\gamma_{1}^{k} \equiv_{n+1} [s_{1}^{k}, \ldots, s_{n}]$.
\noindent
Further, note that there exist elements $\sigma_{1}, \ldots, \sigma_{m} \in \mathcal{S}$ with $m < N_{0}$ and an element $\gamma \in \Gamma_{n}$ such that  $$s_{1}^{k} = \sigma_{1} \ldots \sigma_{m} \gamma$$
The elements $\sigma_{1}, \ldots, \sigma_{m}, \gamma$ depend on $\gamma_{1}$ and $k$, but their number does not. Now, we have: $$\gamma_{1}^{k} \equiv_{n+1} [\sigma_{1} \ldots \sigma_{m} \gamma, \ldots, s_{n}]  \equiv_{n+1} [\sigma_{1} \ldots \sigma_{m}, \ldots, s_{n}],$$
\noindent
 where the last equality stems from the centrality of $\Gamma_{n}$. The last expression is a word in the elements of $\mathcal{S}$, whose length is bounded from above by a number that does not depend on $k$. This fact is true not just for $\gamma_{1}$, but for $\gamma_{2}, \ldots, \gamma_{p}$. Since the group $A_{n}$ is abelian, and every element in it can be written as a product of powers of $\gamma_{1}, \ldots, \gamma_{p}$, we get that $A_{n}$ is generated by $\mathcal{S}$ with finite diameter, as required. $\Box$

\paragraph{Proof of Theorem \ref{thm1.2}.} Theorem \ref{thm1.2} is a direct consequence of Lemma \ref{BoundedGenerationLemma} and induction. $\Box$
\paragraph{Proof of Theorem \ref{thm1.1}.} Let $H = H_{1}(S,\mathbb{Z})$. There exists a simple closed curve in $\pi$ that is mapped to $e_{1}$ under this mapping. The mapping class group acts on $H$, and it is well known that this action induces a surjective homomorphism onto $Sp_{2g}(\mathbb{Z})$ (\cite{primer}). Furthermore, the non-separating simple closed curves form a single mapping class group orbit. Thus, by Lemma \ref{SympLemma} and Theorem \ref{thm1.2}, $\pi$ is nilpotent-bounded with respect to $\mathcal{S}$. $\Box$
\paragraph{Proof of Corollary \ref{cor1.3}} This is a direct result of Theorem \ref{thm1.2}, and the fact that any element of $x \in H$ can be written as $x = h_{1} + \ldots + h_{k}$ with $h_{i} \in H_{i}$.

\section{Further Questions.}
The contrast between the result in this paper and Calegari's result that $\pi$ has infinite diameter with respect to $\mathcal{S}$ gives rise to several questions.

\begin{q1}
Recall that $L_{n} = \pi / \pi_{n}$. By Theorem \ref{thm1.1}, $L_{n}$ has finite diameter with respect to $\mathcal{S}$. Call this diameter $d_{n}$. The sequence $\{d_{n} \}_{n = 1}^{\infty}$ is nondecreasing. Is this sequence bounded? If so, by what value. If not, what is its asymptotic growth rate?
\end{q1}

\begin{q2}
The lower central series is but one of the important series of nested subgroups of $\pi$. Another such series is the derived series, whose elements are quotiens of surjections onto solvable groups. This sequence is defined by:
$$\Gamma^{(0)} = \Gamma, \Gamma^{(n+1)} = [\Gamma^{(n)}, \Gamma^{(n)}] $$
\noindent Is the conclusion of Theorem \ref{thm1.1} if we replace the word nilpotent with the word solvable?
\end{q2}



\noindent
Dept. of Mathematics, University of Chicago \\
5734 University Ave. Chicago, IL 60637 \\
E-mails: khalid@math.uchicago.edu, asaf@math.uchicago.edu

\end{document}